\documentclass[times]{article}
\usepackage{natbib}
\usepackage{lineno,hyperref,enumitem,soul}
\usepackage{mathrsfs}
\usepackage{booktabs}
\usepackage{amsmath}
\usepackage{algorithm}
\usepackage[noend]{algpseudocode}
\usepackage{placeins}
\usepackage[pdftex,dvipsnames]{xcolor}  
\usepackage{subcaption} 
\usepackage[colorinlistoftodos,prependcaption]{todonotes}

\usepackage{tikz}
\usetikzlibrary{calc,trees,positioning,arrows,chains,shapes.geometric,%
    decorations.pathreplacing,decorations.pathmorphing,shapes,%
    matrix,shapes.symbols}

\tikzset{
>=stealth',
  punktchain/.style={
    rectangle, 
    rounded corners, 
    draw=black, very thick,
    text width=10em, 
    minimum height=3em, 
    text centered, 
    on chain},
  line/.style={draw, thick, <-},
  element/.style={
    tape,
    top color=white,
    bottom color=blue!50!black!60!,
    minimum width=8em,
    draw=blue!40!black!90, very thick,
    text width=10em, 
    minimum height=3.5em, 
    text centered, 
    on chain},
  every join/.style={->, thick,shorten >=1pt},
  decoration={brace},
  tuborg/.style={decorate},
  tubnode/.style={midway, right=2pt},
}

\modulolinenumbers[5]

\title{Uncertainty in marine weather routing}

\author{Thomas Dickson \and Helen Farr \and David Sear \and James I R Blake}

\begin{document}



\maketitle





\section{Introduction}

\begin{abstract}
    Weather routing methods are essential for planning routes for commercial shipping and recreational craft. This paper provides a methodology for quantifying the significance of numerical error and performance model uncertainty on the predictions returned from a weather routing algorithm. The numerical error of the routing algorithm is estimated by solving the optimum path over different discretizations of the environment. The uncertainty associated with the performance model is linearly varied in order to quantify its significance. The methodology is applied to a sailing craft routing problem: the prediction of the voyaging time for an ethnographic voyaging canoe across long distance voyages in Polynesia. We find that the average numerical error is $0.396\%$, corresponding to $1.05$ hours for an average voyage length of $266.40$ hours. An uncertainty level of $2.5 \%$ in the performance model is seen to correspond to a standard deviation of $\pm 2.41-3.08\%$  of the voyaging time. These results illustrate the significance of considering the influence of numerical error and performance uncertainty when performing a weather routing study.
\end{abstract}

Marine weather route modelling identifies the minimum time path between two locations given the performance model and the weather that is expected on a given route. Route modelling can be used within the design process for new designs or to provide routes for existing marine craft. When applied operationally, weather routing is essential for improving the safety of mariners at sea through identifying routes which minimise fuel costs, risk of harm or time.

Quantifying the significance of uncertainty in the weather routing process allows the credibility of the predictions to be quantified, managing the expectations of the operators with regards to the accuracy of the supplied route. Key sources of uncertainty are the weather data, the accuracy of the performance model and the numerical error in the solution of the shortest path algorithm. However, research has not investigated the quantification of the numerical error of the routing algorithm or the impact of uncertainty in the performance model. The increasing availability of high performance computing and associated improvement of programming languages allows previously computationally intractable problems to be simulated. 

This paper introduces a methodology for quantifying the significance of numerical error and performance uncertainty in a marine weather routing study. Initially the numerical error of the shortest path algorithm will be quantified through varying the discretization of the environment. The influence of uncertainty of the performance is estimated through varying the performance linearly about the original performance model.

A marine weather routing problem with a high level of uncertainty is the modelling of ethnographic voyaging canoes completing colonisation voyages across Polynesia. This problem involves predicting the voyaging time of a sailing craft over a given route for a range of weather conditions, a typical design problem. It is possible to model the performance of any marine craft as a function of wind and wave conditions, allowed this method to be applied to any marine craft.

\subsection{Literature review}
\label{ssec:route_modelling}

This review will discuss the literature on sailing craft performance prediction and weather routing algorithms. The weather routing process predicts the optimum route for a sailing craft to take between two points. The method considers the performance model of the sailing craft, the shortest path algorithm and the environmental data used to identify the optimum path. 

The performance of a sailing craft is determined by the balance of the driving force generated by the wind passing over the sail against the resistance of the hull and appendages. The wind passing over the sail also generates a heeling moment which is balanced against the righting moment of the hull. The interaction of these key force balances and the additional moments determine the speed and heel angle of the sailing yacht \citep{Philpott1993}.

Static and dynamic velocity prediction programs (VPPs) are used to model the performance of a sailing craft. Static VPPs predict the conditions at which the forces and moments acting on a sailing craft are balanced \citep{Philpott1993}. Dynamic VPPs evaluate the forces acting on the sailing craft over a series of time steps, within the context of a short race \citep{Philpott2004}. Static VPPs require less information than dynamic VPPs on the design of a sailing craft to provide performance predictions but are considered to be less accurate as a result.

The first research into solving the sailing craft route optimisation problem used a recursive dynamic programming formulation which divided the domain into nodes over which the shortest path was calculated \cite{Allsopp1998a}. Different wind models \cite{Philpott2004,Dalang2015} or race strategy and opponent models \cite{Spenkuch2014,Tagliaferri2017} have been used to improve the accuracy of sailing routing models. 

The influence of different methods of modelling the ability for a sailing craft to sail upwind has been explored \cite{Stelzer2008} along with modelling the time taken to complete course changes \cite{Ladany2017}. 
Sailing craft race modelling has typically minimised either the time taken to complete a course \cite{Ferretti2018}, risk of losing to an opponent \cite{tagliaferri2014,Spenkuch2014} or reliability \citep{Dickson2018}. However, the majority of sailing craft routing research only applies to a short course sailed over a small spatial and temporal domain.

The long course modelling problem can be characterised through its consideration of larger spatial and temporal domains over which the shortest path is solved. Typical drawbacks of short course routing methods involve their requirement for a predictable wind field and the requirement for the entire course to be modelled as being flat with respect to the curvature of the earth. The haversine formula becomes accurate at distances longer than $1$ nm \citep{Allsopp1998}, this provides an indication at the crossover between short and long course routing methods.

The marine weather routing literature has developed a range of different shortest path algorithms. The different approaches were classified into calculus of variations, grid based approaches and evolutionary optimisation \citep{Walther2016}. A time dependent approach based on the calculus of variations solved the shortest path problem through identifying the shortest path to take through considering a series of fronts reached within incremental time steps \citep{Bijlsma1975}. 
Forward dynamic programming approaches involve recursively solving the shortest path over a grid of locations generated along the great circle line between the start and finish locations \citep{Allsopp1998}. This approach can be two dimensional or three dimensional considering how the variables such as time and fuel cost are optimised \cite{hagiwara1987}. A claimed improvement is to use a floating grid system which updates potential locations every time step \citep{Fang2013}. These grids can be generated based off local wind conditions \citep{Tagliaferri2017}, although this approach would be challenging for larger spatio-temporal domains as found in long course routing.

Another route modelling approach has involved generating candidate routes using biased Rapidly exploring Random Trees which solve for the minimum energy path using the A$^*$ algorithm \citep{Rao2009}. 
A multi objective genetic algorithm has been applied to solve for the optimum path considering multiple safety and fuel constraints \citep{Hinnenthal2008}. An improvement came with the implementation of a fuzzy logic model to model the performance of a sail driven vessel \citep{Marie2014}.

A method of iteratively aggregating the shortest path over multiple weather scenarios has been introduced \citep{Allsopp1998}. The approach of using ensemble weather scenarios has been shown to be more accurate than single scenarios with application to marine weather routing \citep{Skoglund2015}. Multiple algorithms and objective functions have been implemented and applied to the marine weather routing problem. To the authors knowledge there has been no study of how uncertainties within the solution algorithm and objective function influence the confidence that may be had in the final result.

\subsection{Uncertainty analysis}
\label{ssec:uncertainty_analysis}

Uncertainty analysis is the process of quantifing how likely certain states of a system are given a lack of knowledge on how certain parts of the system operate. The choice of how to categorize and thus simulate uncertainty is significant in the scientific modelling process \citep{Kiureghian2007}. In the marine environment it is possible to classify uncertainty as being aleatory uncertainty or epistemic uncertainty \citep{Bitner-Gregersen2014}. Aleatory uncertainty is the inherent randomness in a particular parameter, it is not possible to reduce this. Epistemic uncertainty is knowledge based, it is possible to reduce this quantity through collecting more information on the process in question. Epistemic uncertainty can be broken down into data uncertainty, statistical uncertainty, model uncertainty and climate uncertainty. Data uncertainty is associated with the error associated with collecting data from experiment or model. Statistical uncertainty is a consequence of not obtaining enough data to model a given phenomenon. Climate uncertainty addresses the ability for the climate variables over a given spatial-temporal domain to be representative given the nature of the weather and climate change \citep{Bitner-Gregersen1990}.

The weather has been described as a chaotic process \citep{Lorenz1963}. The chaotic nature of the weather limits the ability to utilise numerical models to predict into the future or to model what occured in the past. A chaotic system can be thought of as one where the present state determines the future state, but the approximate present state doesn't determine the approximate future. The use of ensemble weather scenarios generated using different intial conditions attempts to mitigate the associated uncertainty with weather forecasts and reanalysis data \citep{Slingo2011}. Given that the weather is a chaotic process, it is likely that any solution process for the shortest path sequential decision making process may identify one or more stable solutions. Currently the simplest method of simulating weather uncertainty is through using past weather data.

The error in a scientific model is the accuracy at which it estimates a real system. The key source of error is the ability for the shortest path algorithm to identify the optimal path based on the discretization of the environment. Typical discretization error calculation methods require grids of significantly different sizes to be solved for a given set of initial conditions \citep{Roache1997}. Through examining the rate of convergence of solutions from different grids it is possible to extrapolate the solution for an infinite discretization.

\FloatBarrier

\section{Method}
\label{sec:method}

\begin{figure}
    \centering
    \begin{tikzpicture}
        [node distance=.4cm,
        start chain=going below,]
           \node[punktchain, join] (intro) {Initial conditions} ;
           \node[punktchain, join] (ra) {Routing algorithm};
           \node[punktchain, join] (Numerr) {Numerical error};
           \node[punktchain, join] (perf) {Performance uncertainty};
            \node[punktchain, join] (uncq) {Uncertainty quantification};
    \end{tikzpicture}
    \caption{Method used to quantify uncertainty in marine weather routing.}
    \label{fig:method_figure}
\end{figure}
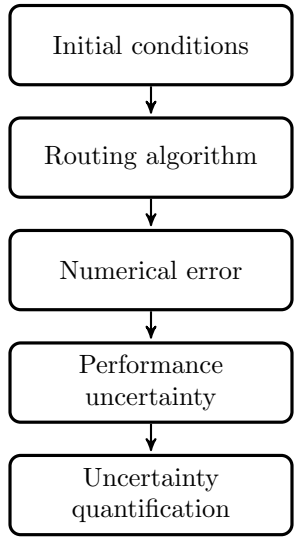

Figure \ref{fig:method_figure} shows the method used in this research. The initial conditions such as the route and environmental conditions are specified. A suitable routing algorithm is identified and implemented. The numerical error is then estimated and the influence of performance uncertainty is modelled. The method is concluded with an analysis of the impact of performance uncertainty and numerical error on the interpretation of the final set of voyaging time results.

\subsection{Routing model}
\label{ssec:routing_algorithm} 

This paper uses the Dynamic Programming (DP) algorithm as applied to long distance sailing craft in \citep{Allsopp1998}, the simplicity of this algorithm means that the source of uncertainty lies in the discretization of the domain over which the shortest path is solved. 

\begin{figure}[!ht]
    \includegraphics[width=0.9\linewidth]{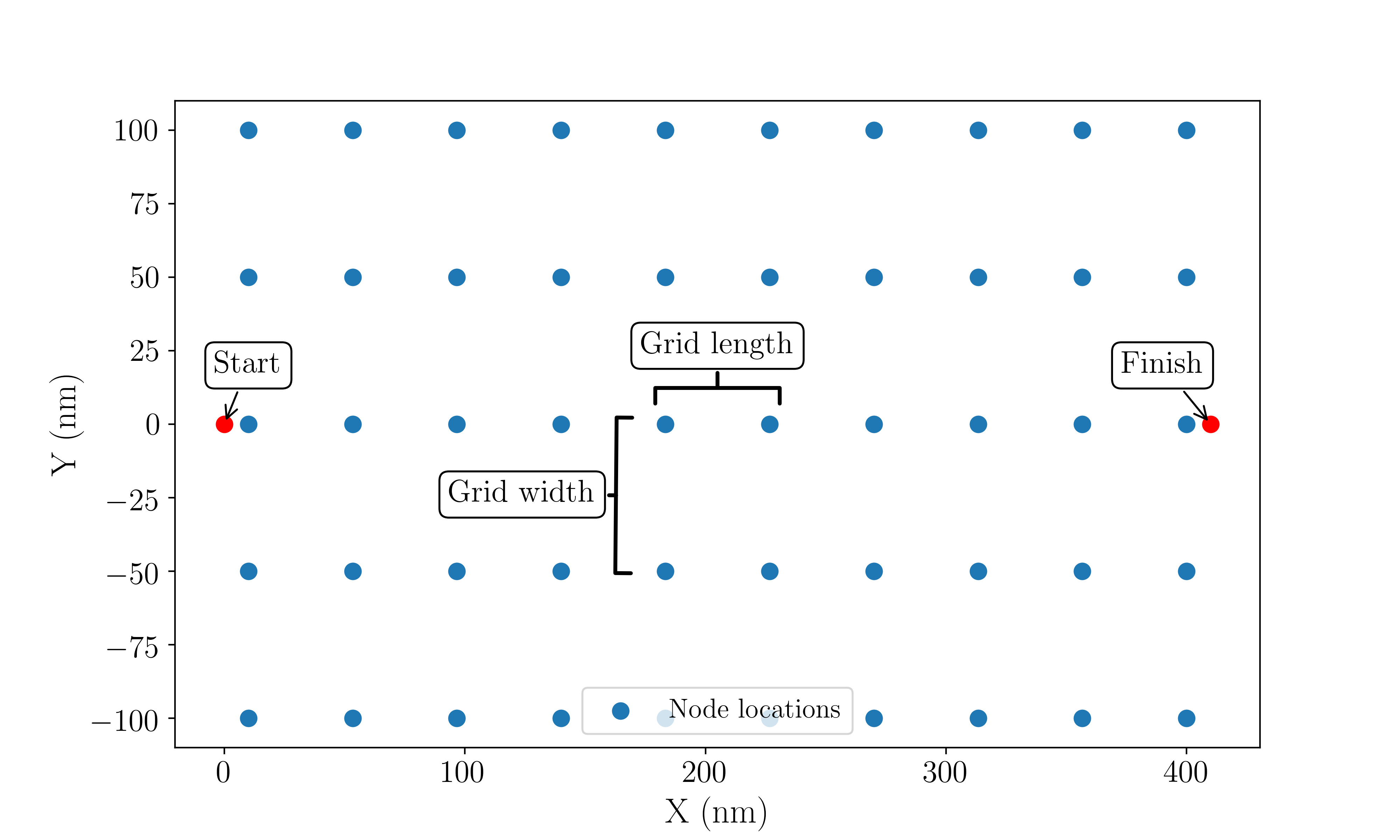}
    \caption{Discretized domain along great circle line between voyage start and finish.}
    \label{fig:discretized_domain}
\end{figure}

We consider a large domain on the Earths surface over which the modelled sailing craft could hypothetically sail. This domain is discretized into an equal number of locations both in parallel and perpendicular to the Great Circle drawn between the start and finish locations. The distance between each location is controlled through defining the maximum distance between nodes, $d_n$. $d_n$ can be controlled seperately as the grid height or grid width, as seen in Figure \ref{fig:discretized_domain}, for this research it is set as being equal. This generates a grid of nodes with equal numbers of ranks as well as nodes within each rank.

For the position at any given node $i$  the travel time between nodes $i$ and a node on the next rank $j$ along the arc $(i, j)$ starting at time $t$ is $c_{arc}(i, j, t)$. The cost function, $c_{arc}(i, j, t)$, provides an estimate of the time taken for a sailing craft to sail between two points given the environmental conditions at time $t$. An intial speed estimate is taken from interpolating the results of a performance prediction analysis for the specific wind condition. This speed is then modified in order to account for the wave conditions being experienced. A final optimisation takes place to identify the optimum heading to sail between the two points given the current the craft is exposed to. If it is not possible to achieve a positive speed towards location $i$ from location $j$ then the speed will be set to $0$ returning an infinite travel time for that particular arc. It is possible to penalise areas of the domain in this manner and identify combinations of initial conditions which are unable to return 

The minimum time path is identified using a forward looking recursive algorithm, described in Equation \ref{eq:forward_recursion}. $f^*(i, t)$ is the time taken for the optimal sequence of decisions from the node-time pair $(i, t)$ to the finish node and $j^*(i, t)$ is the successor of $i$ on the optimal path when in state $t$. Through solving $j^*(i, t)$ within $f^*(i, t)$ it's possible to solve for the shortest path. $\Gamma_i$ is the set of all nodes on the graph. The algorithm iterates from the start node to the finish node and updates each node in between with the earliest time that it can be reached.

\small
\begin{equation} \label{eq:forward_recursion}
    \begin{split}
        f^*(i, t) = \begin{cases}
            0, & i = n_{finish} \\
            \min_{j \in \Gamma_i} [ c_{arc}(i, j, t) + f^*(j, t+c_{arc}(i, j, t))], & \text{otherwise}
    \end{cases} \\
    j^*(i, j) = \text{arg} \min_{j \in \Gamma_i} [c_{arc}(i, j, t) + f^*(j, t+c_{arc}(i, j, t))], i \neq n_{finish}
    \end{split}
\end{equation}
\normalsize

\subsection{Uncertainty simulation}
\label{ssec:unc_sim}

The primary sources of uncertainty in sailing craft route modelling lie in the accuracy of the solution algorithm, the performance model of the sailing craft and the weather. This section introduces the method used to investigate the solution accuracy of the algorithm and the influence of performance uncertainty of a sailing craft. 


The discretization of the environment is the primary source of numerical uncertainty in the solution algorithm. The key parameter determining the fidelity of the simulation is the distance between nodes $d_n$. We are interested in the significance of this parameter on the results on any given routing analysis. In order to consider the influence of the discretization of the domain on sailing craft routing results the lower limit of $d_n$ must be identified which will allow two coarser grid sizes to be chosen.

The lower limit of $d_n$ is selected based on the minimum distance at which the cost function retains accuracy.
Influencing factors include the discretization of the sailing domain and the spatial and temporal resolution of the weather data used. For example, if the time taken to travel between two locations is greater than the temporal resolution of the weather dataset then there are changes in weather conditions which are not being modelled. 

The weather data used is from the ERA20 climate model has a spatial resolution of $125$ km and a temporal resolution of $3$ hours \citep{Poli2016}. The weather data is linearly intepolated to the discretization of the nodes in the environment. If the distance between nodes is much larger than the original spatial discretization of the weather data then the change in the weather conditions will not be fully modelled. The voyaging time is used to look up the closest weather data time step.


The accuracy of the performance model is limited to certain scales, consequently it may not make physical sense to apply such a model below a certain length of time or distance. Despite this, a quantification of the numerical solution algorithm is still required in order to identify whether the shortest path solution is stable. One index used to quantify numerical uncertainty in CFD is the grid convergence index (GCI) \citep{Roache1997,Division2008} which has only recently been applied to the sailing craft routing problem \citep{Dickson2018}. 

The method of numerically quantifying error using the Grid Convergence Index (GCI) is fully described in \cite{Division2008}. It involves the solution of the algorithm over multiple discretizations of the environment which exponentially increase in detail. The grid height, $h$, is the measurement unit of grid size and is calculated using Equation \ref{eq:h_calculation}. $\Delta A_i$ is the size of the $i$th cell and $N$ is the total number of the cells used for computation. The solution trend from the three distinct grid sizes is extrapolated towards $h \rightarrow 0$ where the extrapolated solution is used to estimate the associated discretization error.

\begin{equation}
    h = \Big[ \frac{1}{N} \sum_{i=1}^N (\Delta A_i) \Big]
    \label{eq:h_calculation}
\end{equation} 

We wish to quantify whether uncertainty surrounding the performance modelling of a sailing craft is significant within the context of sailing craft route modelling. To begin the process the performance model is varied linearly about its original performance value. Each generated performance model is known as $P_{unc}$, where $unc$ is the percentage that the original performance is varied by.

Algorithm \ref{alg:unc_sim_routine} describes the uncertainty simulation routine for a specific route at a specific time. It shows how the shortest path is simulated for each $d_n$ parameter and range of $P_{unc}$ values at a start time $t$ for a given route. The number of start times and uncertain performances to be simulated is dependent on the computational resources available. The code for this method was implemented in the programming language Julia \citep{Bezanson2017}. 

\begin{algorithm}[!ht]
    \begin{algorithmic}[1]
    \Procedure{Routing uncertainty simulation}{}
    \For{$d_n$ in $[h_1, h_2, h_3]$}
        \State{Generate discretized environment}
        \State{Spatially interpolate wind, wave and current data for each node}
        \For{$P_{unc} in [50\%, ..., 150\%] \times P_{prediction}$}
            \State{$V_{t, d_n, unc} \leftarrow $ Shortest path($t$, $d_n$, $P_{unc}$)}
        \EndFor
    \EndFor
    \EndProcedure
    \end{algorithmic}
    \caption{Route modelling uncertainty simulation routine for a specific start time $t$ and route.}\label{alg:unc_sim_routine}
\end{algorithm}

\FloatBarrier
\section{Application}
\label{sec:application}

The uncertainty route modelling analysis procedure is applied to quantifying the performance of Polynesian voyaging canoes, an application with previously irreducible levels of uncertainty. Modelling the voyaging time for ethographic voyaging craft to complete specific routes will assist understanding how it was possible for Polynesia to be colonised, one of Pacific archaeology's greatest unanswered questions \citep{Irwin2015}. Of interest is the influence of the ENSO oscillation, a key weather phenomenon in the Pacific, on the voyaging time \citep{Montenegro2016a}. One voyaging route of interest is between Upolu and Moorea. Through simulating the influence of performance uncertainty on the voyaging time for a colonisation voyage it will be possible to quantify the influence of seafaring technology on the rate of the colonisation of Polynesia. At the heart of this problem is the solution of a marine weather routing algorithm with a prior requirement to simulate uncertainty.

The shortest path between Tongatapu and Atiu was estimated for a range of different grid sizes, performances and start times. The results from this series of simulations can be used to estimate the difficulty of making the voyage between these two locations using prehistoric seafaring technology. Twenty different performances were generated linearly varying from $-50\%$ to $+50\%$ of the original performance model. Voyages were started every $72$ hours from the 1st January to the 31st December. $1985$ is used to provide weather data for a medium ENSO, an important condition in the context of the study.

\begin{figure}
    \begin{subfigure}[b]{0.5\textwidth}
      \includegraphics[width=\textwidth]{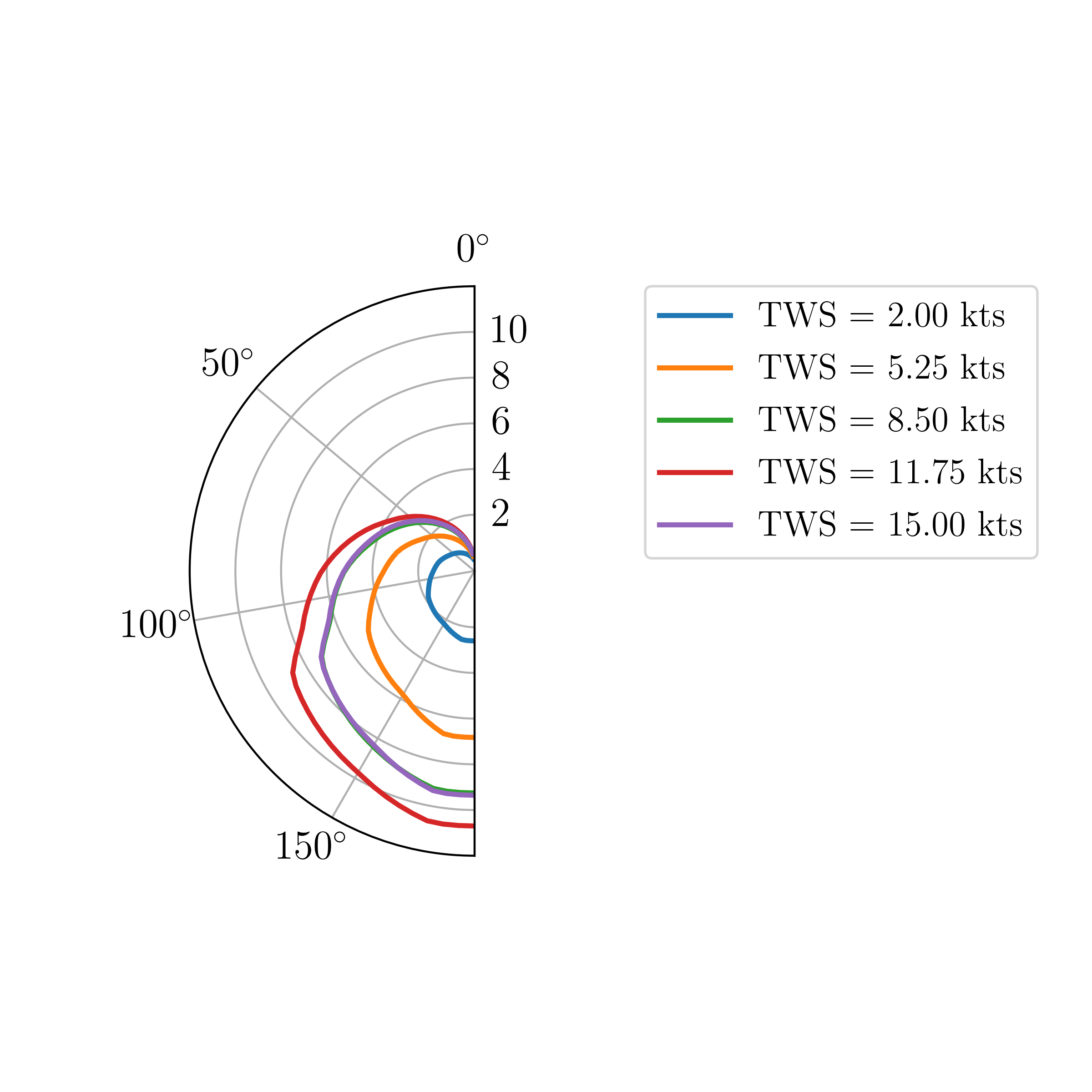}
      \caption{$V_s$ as a function of true wind speed and direction for the ethnographic sailing craft \cite{Boeck2012}.}
      \label{fig:boeck_perf_polars}
    \end{subfigure}
    \begin{subfigure}[b]{0.5\textwidth}
      \includegraphics[width=\textwidth]{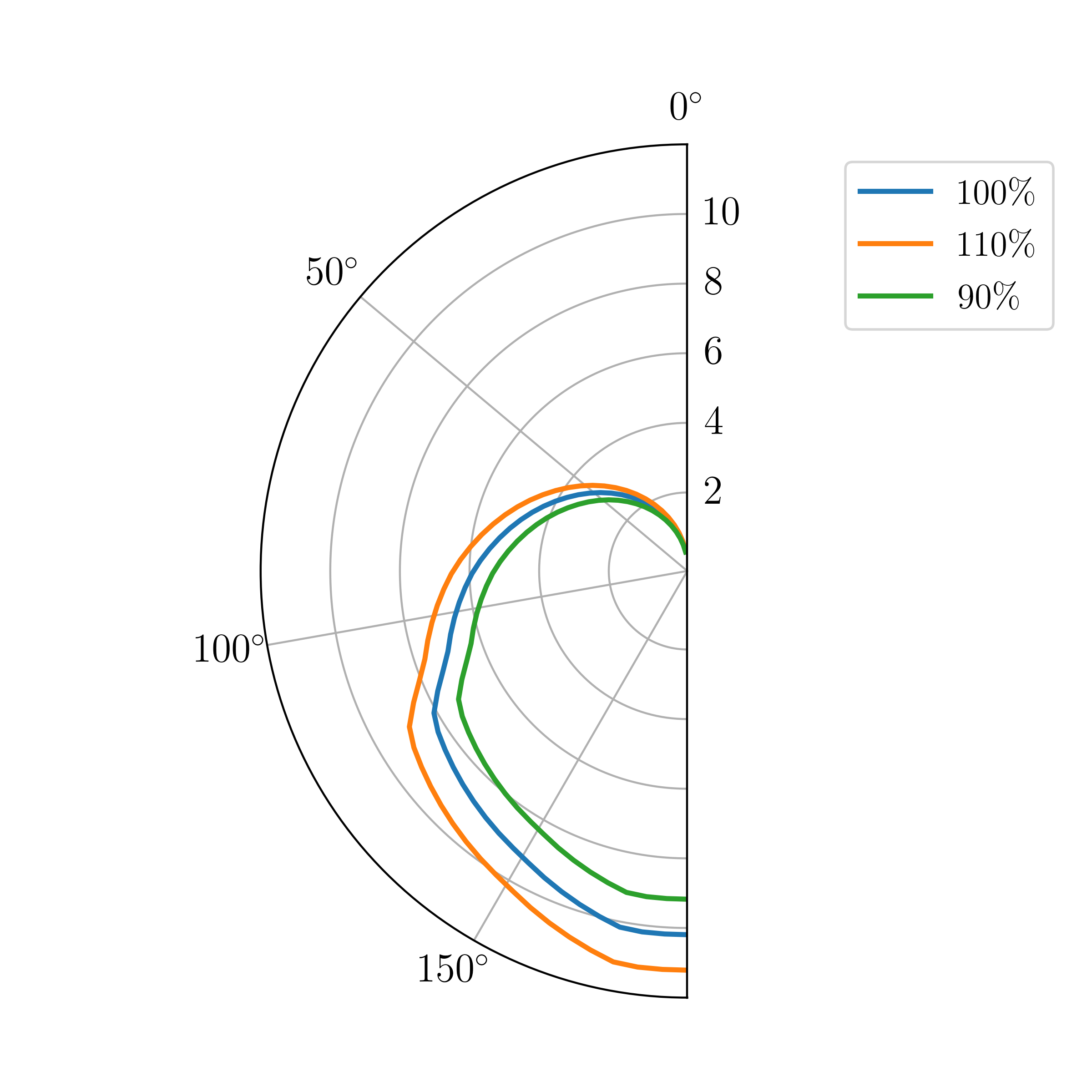}
      \caption{Boat speed, $V_s$, varied $10\%$ about the original performance for a wind speed of $10$ kts.}
      \label{fig:perf_var_example}
    \end{subfigure}
    \caption{Performance of marine vessel used in study.}
\end{figure}

To quantify the efficacy of later voyaging canoe designs, performance data from an outrigger canoe was used \citep{Boeck2012}. The cost function for this craft interpolates the performance from the polar performance diagram, seen in Figure \ref{fig:boeck_perf_polars}. An example of how the performance is varied for a specific wind condition is shown in Figure \ref{fig:perf_var_example}. The wind and wave reanalysis data was downloaded for the year $1982$ from the ECMWF ERA20 model \cite{Poli2016}.  The current data used was sourced from \citep{Bonjean2002}. 

\subsection{Numerical uncertainty}
\label{ssec:results_algorithm_uncertainty}

\subsubsection{Illustration of simulation convergence}
\label{sssec:simulation_convergence}

\begin{figure*}[!ht]
    \centering
    \includegraphics[width=.9\textwidth]{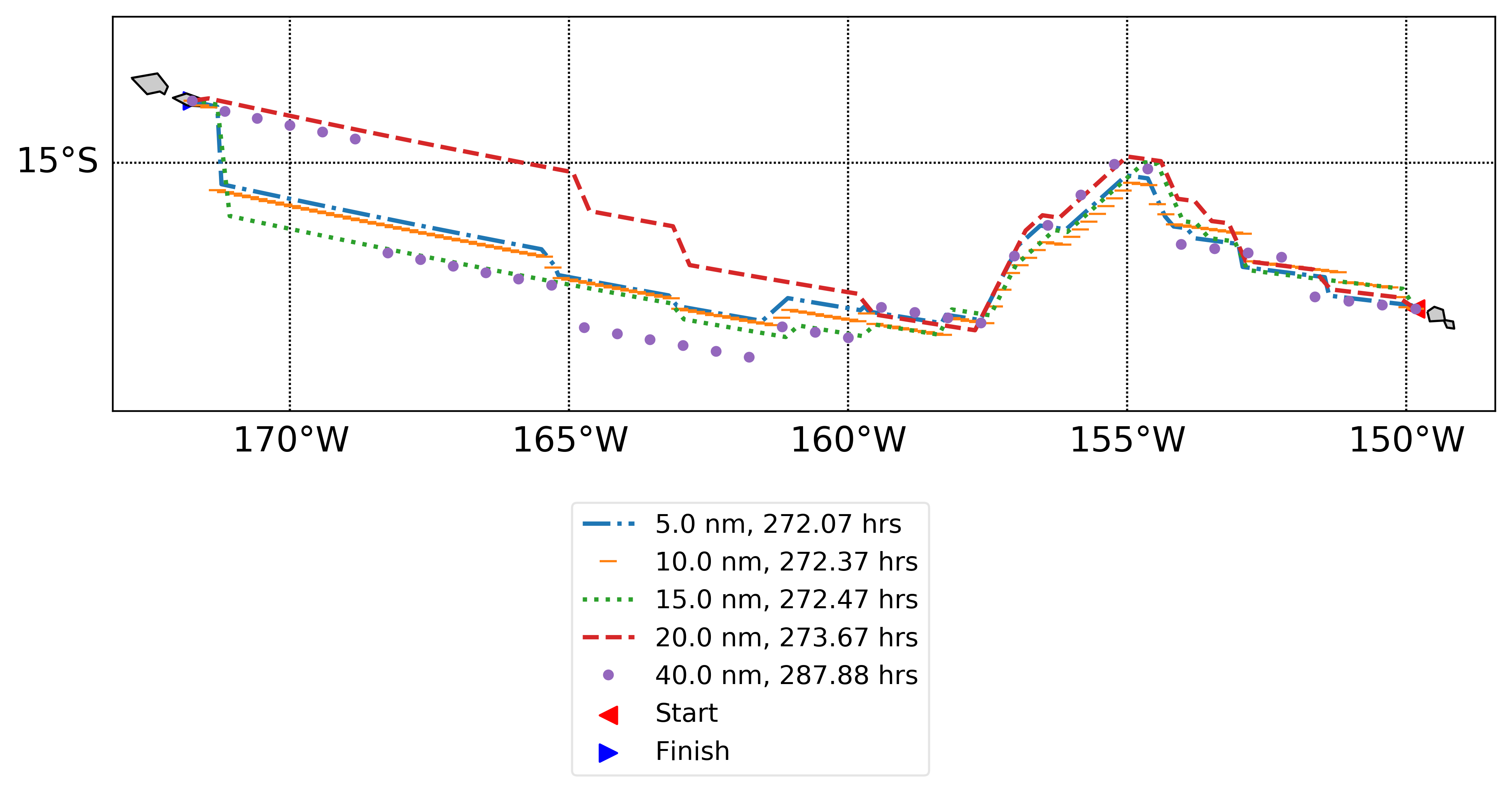}
    \caption{Routes between Upolu and Moorea starting at 00:00 GMT on the 1st January 1985 solved over several different grid widths.}
    \label{fig:grid_width_var}
\end{figure*}

The numerical error in the routing algorithm is a function of the discretization of the domain, parameterised by the grid width, $d_n$. Figure \ref{fig:grid_width_var} illustrates the predicted routes between Tongatapu and Atiu where $d_n$ was reduced in stages between $40.0$ to $5$ nm. It can be seen that as the fidelity of the simulation increases, the voyaging time, $V_t$, reduces. 

\begin{figure}[!ht]
    \centering
    \includegraphics[width=.9\linewidth]{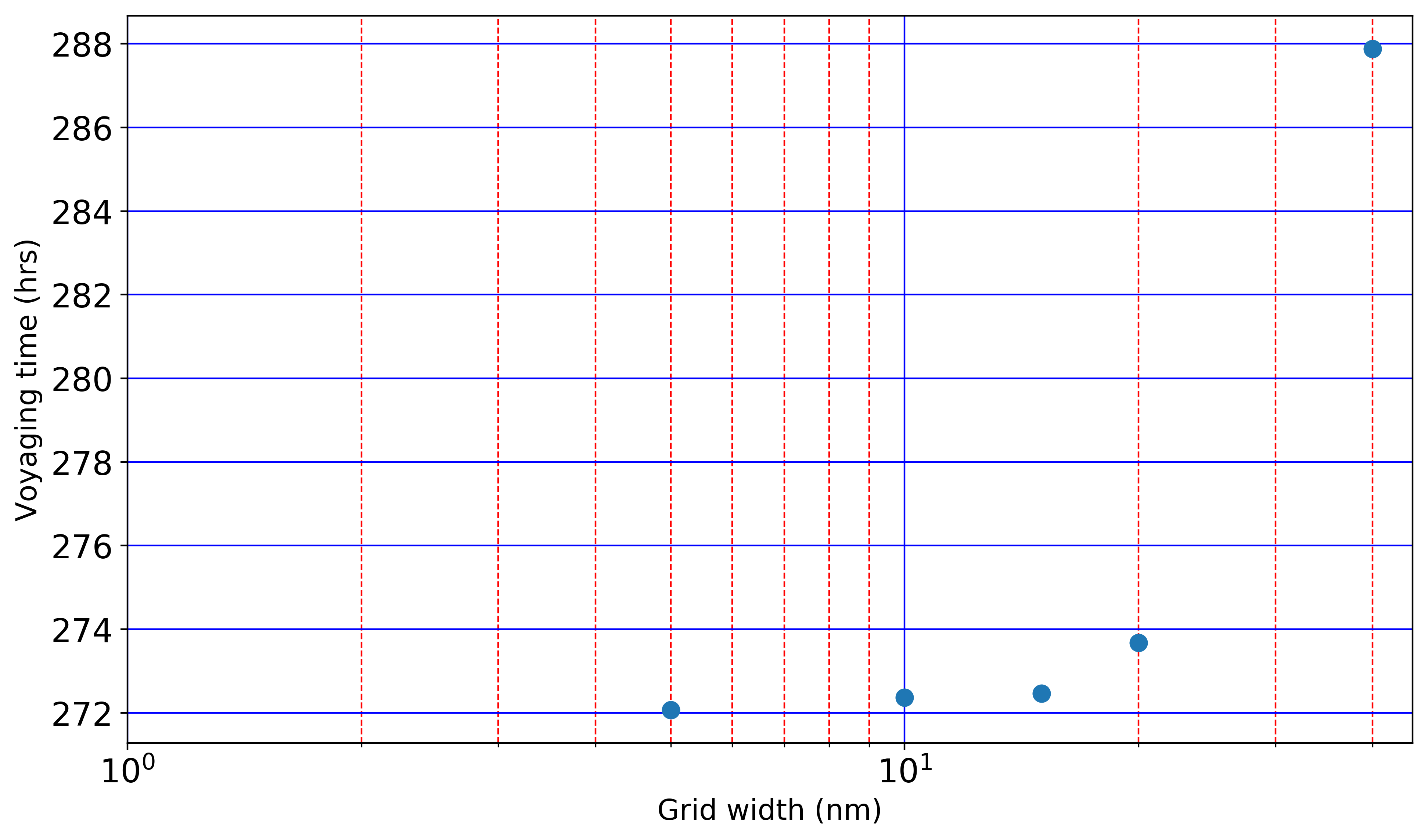}
    \caption{Voyaging time, $V_t$, as a function of grid width for voyages between Upolu and Moorea starting at 00:00 GMT on the 1st January 1985.}
    \label{fig:vt_convergence_plot}
\end{figure}

A convergence plot of $V_t$ as a function of $d_n$ is shown in Figure \ref{fig:vt_convergence_plot}. For the specific initial conditions the relationship is that as $d_n$ reduces, so too does $V_t$. It can be seen that there is a significant change in relationship between the results for $d_n = 40, 20$ and $d_n = 15, 10, 5$ nm and $V_t$. This could be due to the grid discretization becoming significantly finer than the resolution of the original wind and wave data. 

The extrapolated $V_t$, GCI and order of convergence for three combination of heights is included in Table \ref{tab:gci_calcs}. The extrapolated $V_t$ value is where the predicted $V_t$ for the specified heights is extrapolated to $d_n = 0$. For $d_n = 10.0, 15.0, 20.0$ we see there is a low associated GCI as it looks as if the results are converging rapidly. However, we see that for both $d_n = 5.0, 10.0, 15.0$ nm the GCI increases by an order of magnitude as $V_t$ reduces. This indicates that as the accuracy of the simulation improves the application of GCI becomes meaningful. However, there is a large computational cost associated with finer simulations which prohibits the number of simulations required to explore uncertainty in voyaging time. These results show that the GCI calculation method for estimating numerical error must be applied in a manner considering the physical implication of the parameters but also the computational run time of the simulations.

\begin{table}[!ht]
    \centering
    \begin{tabular}{@{}llll@{}}
    \toprule
    \textbf{$d_n$}     & \textbf{$V_t$ (hrs)} & \textbf{GCI} & \textbf{OOC} \\ \midrule
    $5.0, 10.0, 15.0$               & 274.00                                & 0.0014577                     & 2.25                          \\
    $10.0, 15.0, 20.0$              & 272.49                                & 0.0000955                     & 4.34                          \\
    $15.0, 20.0, 40.0$              & 272.64                                & 0.0012429                     & 2.59                          \\ \hline
    \bottomrule
    \end{tabular}
    \caption{Calculation of extrapolated $V_t$, GCI and order of convergence.}
    \label{tab:gci_calcs}
\end{table}

The selection of $d_n$ is a compromise between simulation accuracy and computational run time. Simulation accuracy is determined by the smallest and the largest scales that the cost function can be applied to. The smallest scale is determined by the haversine formula which has significant levels of error for distances below $1$ nm. The largest scale is determined by the rate at which spatio-temporal environmental data becomes available. Parallel computing provides significantly more resources than available previously allowing an increase in the number of different initial conditions that are required to be simulated. However, there are still large costs associated with the finer simulations at $d_n=2.5, 5$.

Figure \ref{fig:vt_convergence_plot} illustrates that as $d_n$ is reduced significantly below the spatial distribution of weather data there is a step change in the solution. The $d_n$ also influences the rate at which new weather data is retrieved and processed. If the journey time for a particular arc between two nodes lasts longer than the time between weather conditions being updated then the solution is being solved over incomplete information. The desire for accuracy must be balanced against computational limitations. There is a cubic relation between the computational run time and fidelity of the simulation. From this set of initial simulations it may be proposed that the $d_n = 10.0, 15.0, 20.0$ nm provide a collection of heights which balance the requirement for accuracy against computational run time.

\subsubsection{Numerical error of voyaging simulations}
\label{sssec:general_numerical_error}

Of interest is the average numerical error for the group of simulations being run. The numerical error must be calculated for each set of initial conditions, but knowledge of the average error will give a sense of the degree of confidence we might have. We are interested in what is the most likely amount of numerical error which might exist.

The order of convergence, GCI value and extrapolated voyaging time were calculated for the times generated as a function of the different grid sizes of $10, 15$ and $20$ nm. The order of convergence measures the rate at which the difference between the magitude of each result changes as a function of the change of the grid width. $95.66\%$ of the remaining results had an order of convergence greater than $1.0$, a necessary requirement in order to extrapolate the result to a grid width of approximately $0$. 

The lack of convergence is due to the solution of a shortest path algorithm, a sequential decision making problem, over a chaotic environment. For those remaining results it is possible to investigate the lack of convergence through running simulations at reduced grid widths. For now, the use of the numerical error reduction procedure helps identifying which combinations of initial conditions result in complex behaviour requiring more analysis.

\begin{figure}[!ht]
    \includegraphics[width=0.9\linewidth]{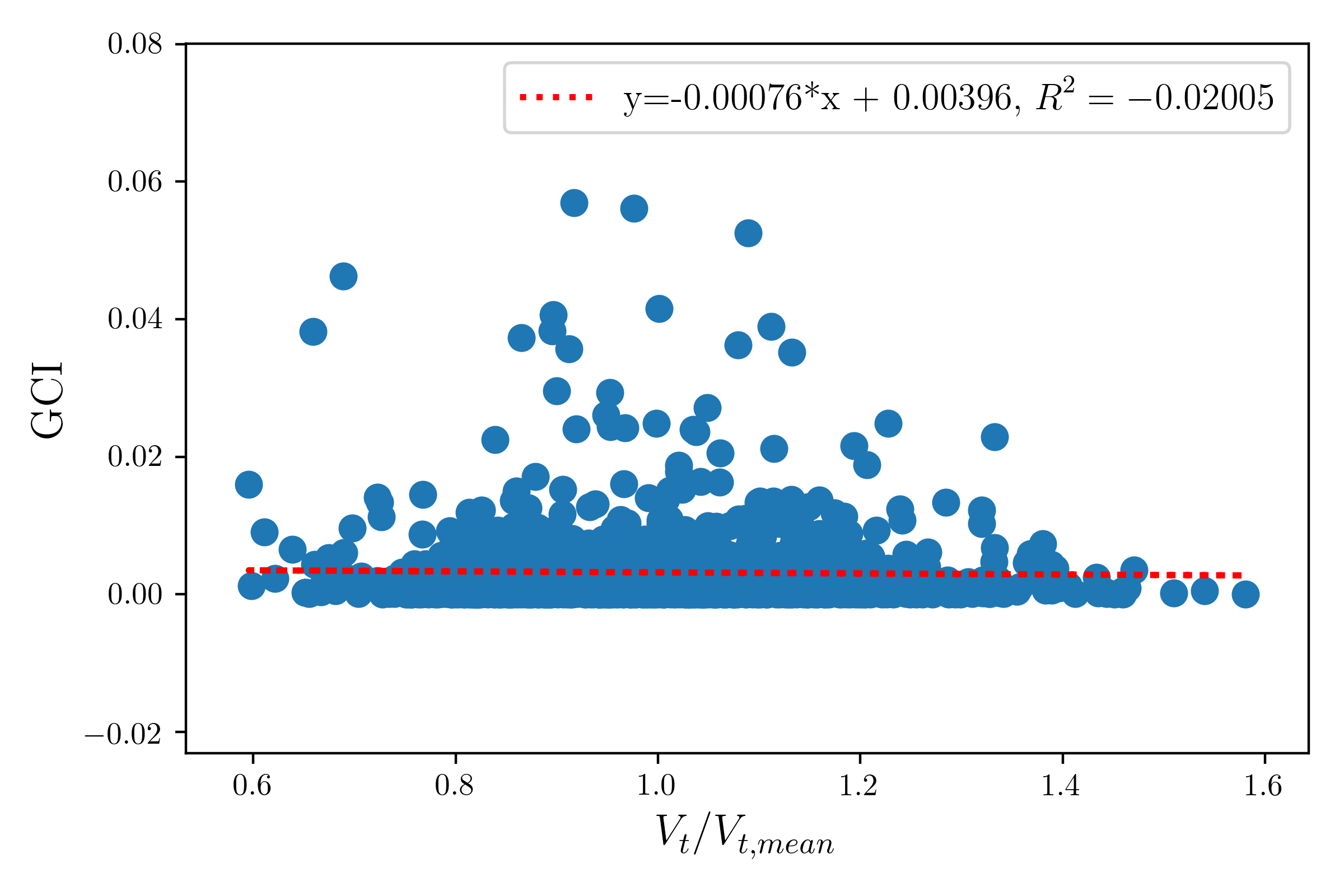}
    \caption{Relationship between GCI index and voyaging time.}
    \label{fig:converged_GCI}
\end{figure}

The relationship between $V_t$ and $GCI$ is shown in Figure \ref{fig:converged_GCI}. It can be seen that there is no relationship between $V_t$ and $GCI$. This can be explained as the same number of calculations are being performed over each grid size. It can be seen that large deviations from the mean $V_t$ are associated with large $GCI$ values. This indicates sets of initial conditions where more analysis is required in order to arrive at credible predictions.
The relationship between the converged results and GCI values estimates an average  of $0.396 \% V_t$ error associated with the whole group of converged voyaging time predictions. This is approximately $1.05$ hrs for the mean voyaging time of $266.40$ hrs. 

\subsection{Performance uncertainty}
\label{ssec:results_performance_uncertainty}

\begin{figure}[!ht]
    \includegraphics[width=0.9\linewidth]{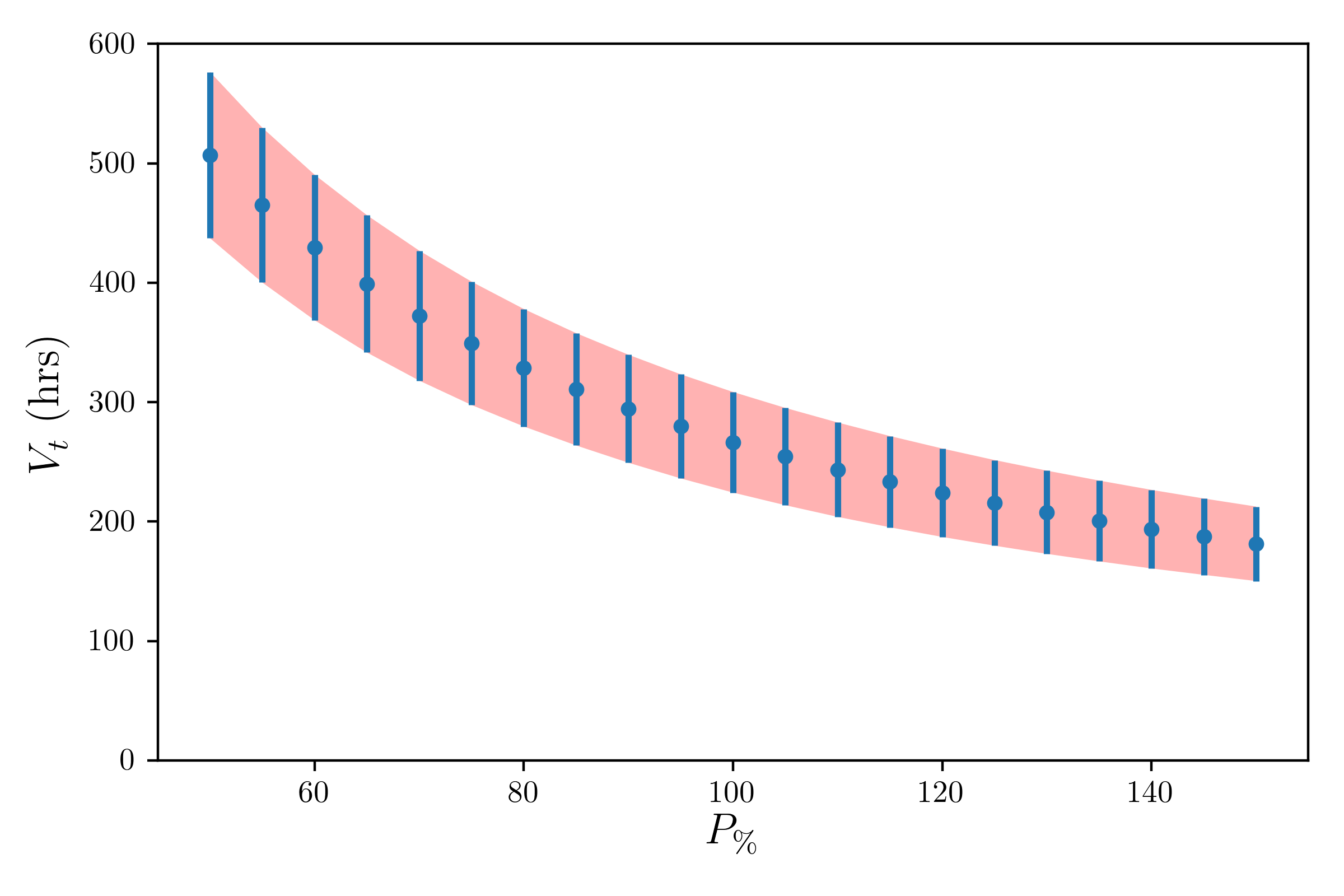}
    \caption{Voyaging time as a function of the different variations of performance. The standard deviation is illustrated using error bars.}
    \label{fig:perf_variation_time}
\end{figure}

The ethnographic voyaging canoe acts as an example of a typical marine craft design problem, albeit one with large levels of uncertainty. Voyages between Tongatapu and Atiu were started every $6$ hours from the $1$st January $1985$ to the 31st December $1985$. The performance model was varied for $21$ steps between $50\%$ and $150\%$ of the original performance. Simulations were solved over a grid size of $d_n = 10$ nm. Significant variation in voyaging time can be seen for any given performance and across all changes in performance, as seen in Figure \ref{fig:perf_variation_time}.

The mean voyaging time appears to be between $230-295$ hours for the unaltered performance with the variation being solely due to the range of environmental conditions. This illustrates that even if there was a high level of confidence associated with the accuracy of the performance model the weather conditions contribute signficantly to the variation of voyaging time.

Across the range of performance variation variation there is a large change from a minimum voyaging time of $98.42$ hours up to a maximum of $682.77$ hours. A lack of confidence in the accuracy of the performance model increases the variation in voyaging time. The accuracy of the performance model must be quantified before performing a routing study so it is possible to quantify the realistic variation in voyaging time.

Of interest is the variance of the voyaging time as a function of the performance variation. Understanding how as the varying of performance influences the voyaging time indicates the degree of confidence that should be held in a given voyaging result, given the confidence in the performance mode. 

\begin{figure}[!ht]
    \includegraphics[width=0.9\linewidth]{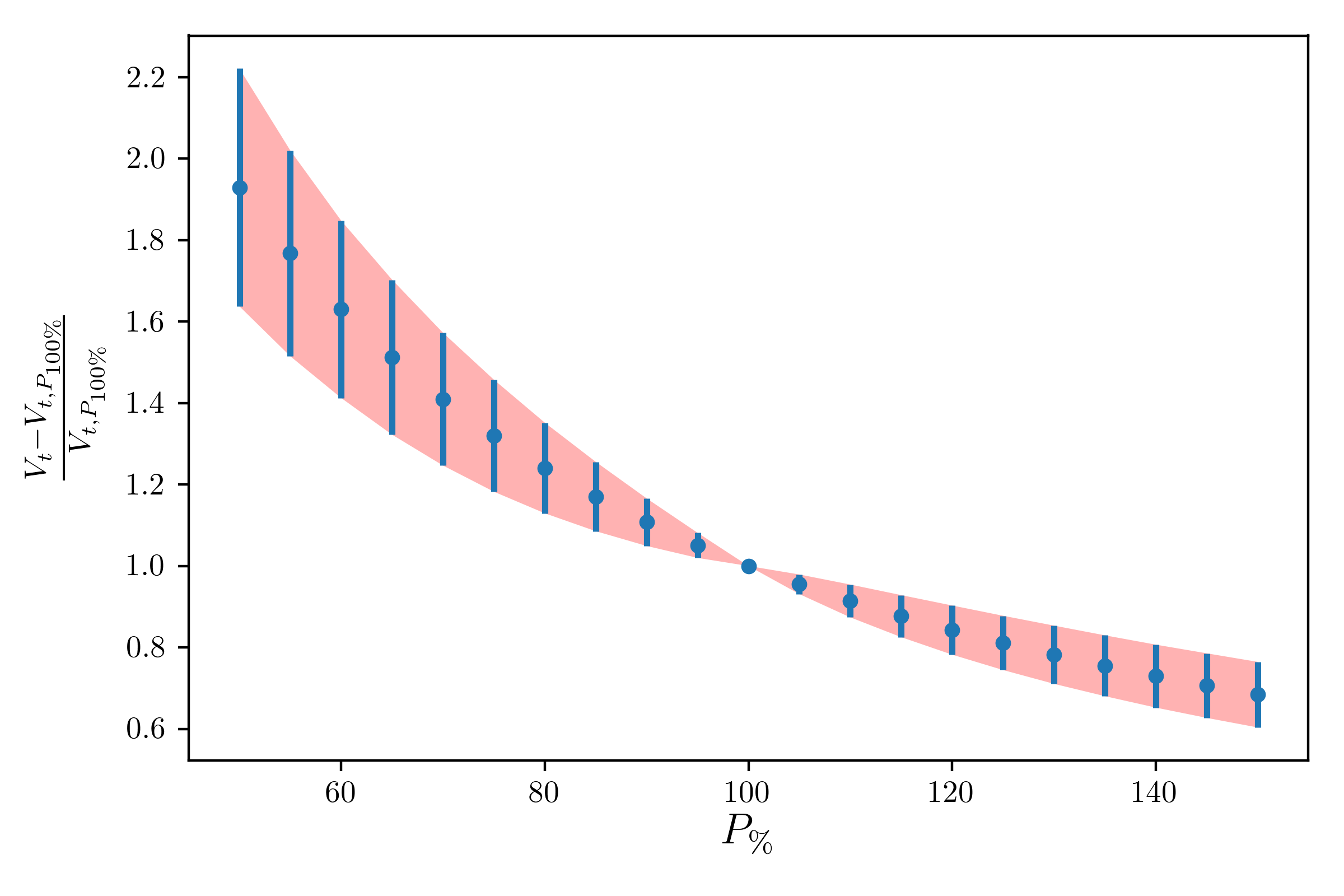}
    \caption{Relationship between the voyaging time and change in performance. The standard deviation is plotted as the error bars.}
    \label{fig:performance_variation}
\end{figure}

Figure \ref{fig:performance_variation} shows the relationship between performance variation and the voyaging time non-dimensionalised with respect to the original performance voyaging time for each start date. This illustrates how variations in performance from the original performance significantly alter the expected voyaging time. It can be seen that as the performance varies from the initial performance the standard deviation of the voyaging time results increases. One key result is that the magnitude of the reduction of voyaging time in response to performance \textit{improvement} is smaller than the magnitude of the increase of voyaging time in response to performance \textit{reduction}. 

There is also a difference between the change in the standard deviation of $V_t$ for equivalent magnitude variations about the original performance. As the performance decreases we see more significant reductions in speed and much larger increases in standard deviation. As the performance increases the magnitude between successive improvements decreases along with a reduction in standard deviation. It can be seen that reductions in performance have more significant negative impacts on voyaging time than equivalent improvements. This non-linear response is due to the slower craft spending more time at sea and consequently being exposed to more variation in the weather.

\begin{figure}[!ht]
    \includegraphics[width=0.9\linewidth]{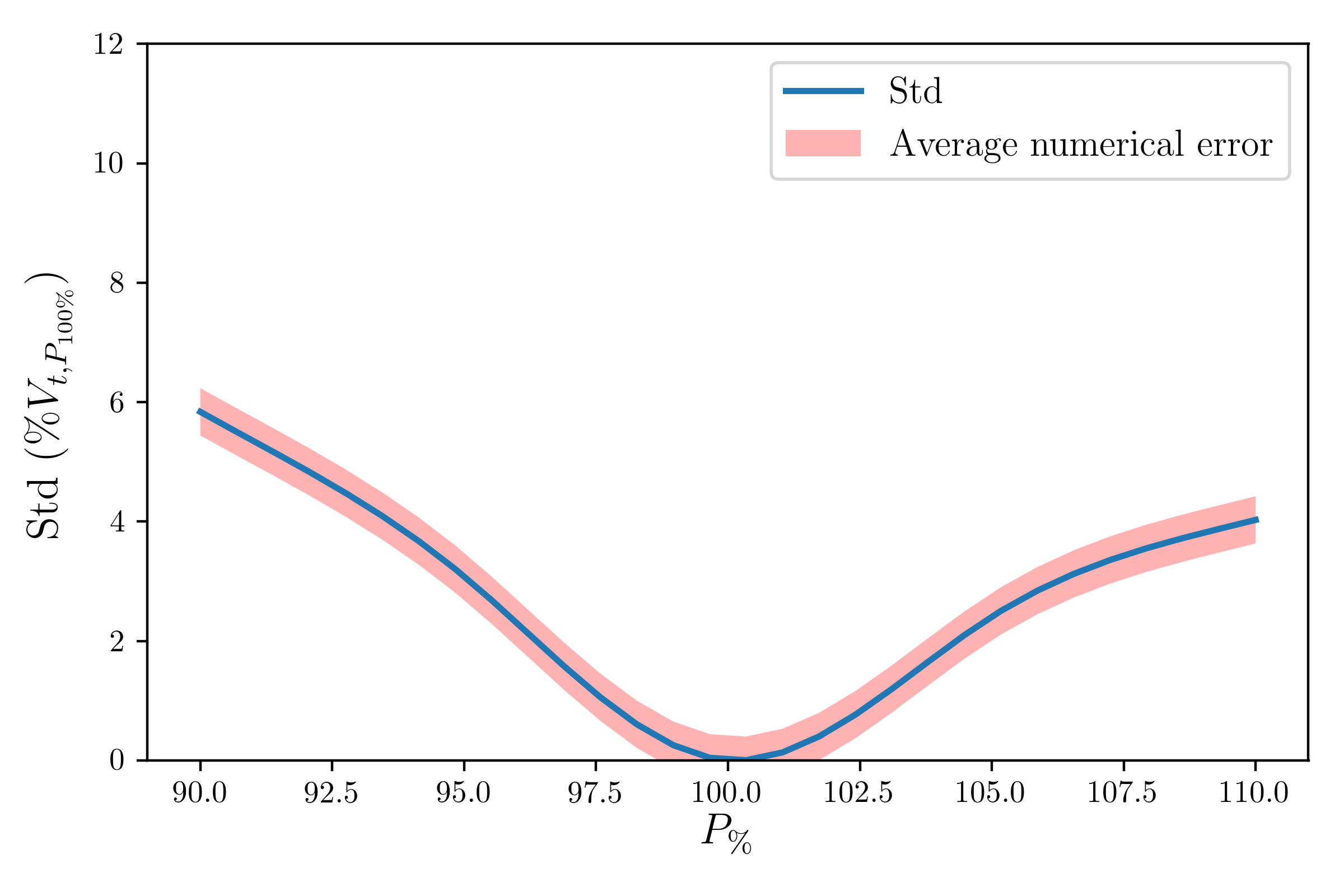}
    \caption{Relationship between performance variation and the standard deviation of $V_{t, P_{100 \%}}$.}
    \label{fig:nd_std_deviation_perf}
\end{figure}

Figure \ref{fig:nd_std_deviation_perf} shows how the standard deviation of the voyaging time varies as a function of the performance variation. The average numerical is overlaid to provide an indication of how significant it may be when using the results of this study. Figure \ref{fig:nd_std_deviation_perf} indicates that the contribution of performance uncertainty is much larger than the numerical error of the algorithm. A variation in performance of $\pm 2.5\%$ causes a standard deviation of $0.8-1.1\%$, equivalent to $2.13 - 2.93$ hrs for the mean voyaging time of $266.40$ hrs. As the variation from the original performance increases to $\pm 5 \%$ we see that the standard deviation increases rapidly to $2.41 - 3.08 \% V_t$, or, $6.42 - 8.21$ hours. These results indicate that uncertainty in the cost function describing the performance of a sailing craft signficiantly change the estimated voyaging time of the shortest path.

The weather conditions are updated every $3$ hours. The magnitude of the numerical error and influence of low levels of performance uncertainty indicate that it is possible for multiple changes in the weather to not be modelled. It would be difficult to quantify the impact of this error on the ability for the weather routing model to approximate the real situation.

\FloatBarrier

\section{Summary and conclusions}
\label{sec:conclusions}

This paper has presented a method for evaluating the impact of numerical error and performance variation in the voyaging time for a marine weather routing problem. This method was applied to a typical design problem; the quantification of the time taken for a sailing craft to complete a specific voyage given uncertainty in its performance. The key results of this study can be summarised as follows;

\begin{enumerate}[noitemsep]
    \item Variation in performance contributed significantly to the variation in voyaging time given the existing variation due to change in environmental conditions.
    \item The numerical error must be calculated for each set of initial conditions. For this problem, $95.3 \%$ of all simulations converged with an average of $0.396\%$ error, equivalent to $1.05$ hours for an average voyage length of $266.40$ hours. 
    \item Slower craft spend more time at sea they are exposed to more variance in the weather conditions, likely contributing towards the non-linear response of voyaging time to performance variation. This means that the uncertainty in the performance model must be quantified to provide credibility to voyaging simulations.
    \item There is a non linear relationship between variation of performance and voyaging time.  The relationship between uncertainty and the standard deviation of voyaging time increases sharply with variations of $2.5 \%$ in performance being associated with standard deviations of $\pm 2.41-3.08\%$ about the mean voyaging time. The influence of uncertainty in the performance model rapidly becomes more influential than the routing algorithm numerical error,
    \item The weather data used updates every $3$ hours. The combination of numerical error and uncertainty in performance model may mean that the approximation of the shortest path is being calculated based off incomplete sets of weather data, or solved using more weather data than would be encountered in practice.
\end{enumerate}

This method of quantifying the numerical error of the solution algorithm and performance uncertainty could be applied to other cases involving marine vessels such as cargo ships. This would allow an understanding of the maximum level of accuracy that could be achieved within commerical practice. Another investigation could be into the uncertainty levels associated with the recorded reanalysis weather data used and how this might influence the result.

Through applying an uncertainty analysis method to the marine weather routing problem we have shown that the influence of performance uncertainty is much larger than any uncertainty associated with the shortest path algorithm used. To provide more accurate routing the uncertainty associated with the performance model used must be reduced. 

\section*{Acknowledgements}

Dr Gabriel Weymouth, Professor Dominic Hudson and Dr Blair Thornton for discussion of methodology and results. Carlos Losada de la Lastra for extensive comments on the text.

\section*{Funding sources}

This work was funded by the Southampton Marine and Maritime Institute and the University of Southampton.


\bibliography{/home/thomas/Documents/PhD/library}

\end{document}